\documentclass{article}

\usepackage{amssymb}

\renewcommand{\thefootnote}{}

\newcommand{\RR}{{\mathbb R}}

\newcommand{\NN}{{\mathbb N}}

\newtheorem{theorem}{Theorem}

\newtheorem{corollary}[theorem]{Corollary}

\newtheorem{lemma}[theorem]{Lemma}

\newtheorem{proposition}[theorem]{Proposition}

\begin{document}

\begin{center}
{\Large \textbf{The natural algorithmic approach of mixed
trigonometric-polynomial problems}}

\bigskip

Tatjana Lutovac${}^{\mbox{\tiny $\,1)$}}$, Branko Male\v sevi\' c${}^{\mbox{\tiny $\,1)\,\ast$}}$,
Cristinel Mortici${}^{\mbox{\tiny $\,2)$}}$
\footnote{$\!\!\!\!\!\!\!\!\!\!\!\!{}^{\ast}\,$Corresponding author,
{\em Telephone}: +381113218321, {\em Fax}: +381113248681}
\footnote{$\!\!\!\!\!\!\!\!\!\!\!$E-mails:
Tatjana Lutovac$\,<${\sl tatjana.lutovac@etf.rs}$>$,
Branko~Male\v sevi\' c$\,<${\sl malesevic@etf.rs}$>$,
Cristinel Mortici$\,<${\sl cristinel.mortici@valahia.ro}$>$}

\bigskip

{\footnotesize \textit{${}^{1)}$Faculty of Electrical Engineering,
University of Belgrade, \\[0.0ex]
Bulevar kralja Aleksandra 73, 11000 Belgrade, Serbia \\[2ex]
${}^{2)}$Valahia University of T\^{a}rgovi\c{s}te, Bd. Unirii 18, 130082 T\^{a}rgovi\c{s}te, Romania; \\[0.2ex]
Academy of Romanian Scientists, Splaiul Independen\c{t}ei 54, 050094 Bucharest, Romania; \\[-0.2ex]
University Politehnica of Bucharest, Splaiul Independen\c{t}ei 313, 060042 Bucharest, Romania}}
\end{center}

\bigskip \noindent {\small \textbf{Abstract.} } The aim of this paper is to
present a new algorithm for proving mixed trigonometric-polynomial
inequalities of the form
\[
\sum\limits_{i=1}^{n}\alpha _{i}x^{p_{i}}\!\cos ^{q_{i}}\!x\sin
^{r_{i}}\!x>0,
\]%
by reducing to polynomial inequalities. Finally, we show the great
applicability of this algorithm and as examples, we use it to
analyze some new rational (Pad$\acute{\mbox{e}}$) approximations of the function $\cos ^{2}{\!x}$,
and to improve a class of inequalities by Z.-H. Yang. The results of our analysis
could be implemented by means of an automated proof assistant, so our work is
a contribution to the library of automatic support tools  for proving various analytic
inequalities.

\bigskip \noindent {\footnotesize \textbf{MSC 2010:} 41A10; 26D05; 68T15; 12L05
41A58}

\bigskip \noindent {\footnotesize \textbf{Keywords:} mixed trigonometric-polynomial functions;
Taylor series; approximations; inequalities; algorithms; automated theorem proving}


\section{Introduction and Motivation}


In this paper, we propose a general computational method for
reducing  some inequalities involving trigonometric functions to the
corresponding polynomial inequalities. Our work has been motivated
by many papers \cite{Zhang_Zhu_2008},
\cite{Klen_Visuri_Vuorinen_2010}, \cite{Mortici_2011},
\cite{Chen_2012},
\cite{Milovanovic_Rassias_2014}, \cite{Z.-H._Yang_2014},
\cite{Bercu_2016} -
\cite{Bercu_2017}
recently published in this area. As an example, we mention
the work of Mortici \cite{Mortici_2011} who extended
Wilker-Cusa-Huygens inequalities, using a method, he called
\emph{the natural approach method}. This method consists in
comparing and replacing $\sin x$ and $\cos x$ by
their  corresponding Taylor polynomials, as follows:%
$$\displaystyle  \begin{array}{ccccc}
\mathop{\mbox{$\sum$}}\limits_{i=0}^{2s+1}%
\displaystyle\frac{(-1)^{i}x^{2i+1}}{(2i+1)!}
\!\!&\!\!<\!\!&\!\! \sin x \!\!&\!\!<\!\!&\!\!
\mathop{\mbox{$\sum$}}\limits_{i=0}^{2s}%
\displaystyle\frac{(-1)^{i}x^{2i+1}}{(2i+1)!}                        \\[2.5 ex]
\displaystyle \mathop{\mbox{$\sum$}}\limits_{i=0}^{2k+1}\frac{(-1)^{i}x^{2i}}{%
(2i)!}
\!\!&\!\!<\!\!&\!\! \cos x \!\!&\!\!<\!\!&\!\!
\displaystyle  \mathop{\mbox{$\sum$}}\limits_{i=0}^{2k} \frac{(-1)^{i}x^{2i}}{%
(2i)!},
\end{array}
$$
for every integers
$s, k \!\in\! \NN_{0}$
and
$x\!\in\!\left(0,\pi/2\right)$.

\break

In this way, complicated trigonometric expressions can be reduced to
polynomial, or rational expressions, which can be, at least
theoretically,  easier  studied (this can be done   using some
softwares for symbolic computation, such as Maple).

\smallskip
For example, Mortici  in  \cite{Mortici_2011}  (Theorem 1), proved the next  inequality:%
\[
\displaystyle
\cos x - \left( \frac{\sin x}{x} \right)^{\!3}
> \,
- \frac{x^{4}}{15}, \;\;\;\;
x \!\in\! \left(0,\mbox{\small $\displaystyle\frac{\pi }{2}$}\right)\!,
\]%
by intercalating  the following  Taylor polynomials, as follows:
\[
\displaystyle
\cos x \,- \left(\frac{\sin x}{x}\right)^{\!3}\!+\,\frac{x^{4}}{15}
>
1 - \frac{x^{2}}{2!}+\frac{x^{4}}{4!}-\frac{x^{6}}{6!} \,-\!
\left(\!\frac{ x - \mbox{\small $\displaystyle\frac{x^{3}}{3!}$} +
\mbox{\small $\displaystyle\frac{x^{5}}{5!}$}}{x} \!\right)^{\!\!3\!}
\! + \, \frac{x^{4}}{15} = \frac{x^{6}R\left(x^{2}\right)}{1728000},
\]%
where  $\displaystyle R\left( t\right) =20000-1560t+60t^{2}-t^{3}.$


\smallskip
\noindent Although  transformation based on the natural approach
method  has been made by several researchers in their isolated
studies, a unified approach has not been given yet. Moreover, it is
interesting to note that just  trigonometric expressions involving
odd powers of  $\cos x$ were studied only, as the natural approach
method cannot be  directly applicable for the function
$\displaystyle \cos^{2}\!x$ in the entire over interval $(0,
\pi/2)$.

\smallskip
The aim of this paper is  to extend and  formalize  the ideas of the
natural approach method  for  a wider class of trigonometric
inequalities, including also those containing even powers of $\cos
x$, with no further restrictions.

 Let $\delta _{1}\!\leq \!0\!\leq \!\delta _{2},$ with $%
\delta _{1}\!<\!\delta _{2}.$ Recall that a function defined by the formula%
\begin{equation}
\label{MTF 1}
f(x)=\sum\limits_{i=1}^{n}\alpha _{i}x^{p_{i}}\!\cos ^{q_{i}}\!x\sin
^{r_{i}}\!x,\ \ \ x\!\in\!(\delta _{1},\delta _{2}),
\end{equation}%
is named a mixed trigonometric-polynomial function, denoted in the
sequel by MTP function \cite{Dong_Yu_Yu_2013}, \cite{Malesevic_Makragic_2015}.
Here, $\alpha _{i}\!\in \!\mathbb{R}\!\setminus \!\{0\}$,
\mbox{$p_{i}, q_{i}, r_{i} \!\in\! \mathbb{N}_{0}$},
\mbox{$n \!\in\! \mathbb{N}$}.
Moreover, an inequality of the form $f(x)\!>\!0$ is called a mixed
trigonometric-polynomial inequality (MTP inequality).

\smallskip
MTP functions currently appear in the monographs on the theory of
analytical inequalities \cite{Mitrinovic_1970},
\cite{Milovanovic_Mitrinovic_Rassias_1994} and
\cite{Milovanovic_Rassias_2014}, while concrete MTP inequalities are
employed in numerous engineering problems (see e.g.
\cite{deAbreu_2009b}, \cite{deAbreu_2012}). A large class of
inequalities arising from different branches of science,  can be
reduced to MTP inequalities. Notwithstanding, the development of
formal methods and procedures for automated generation of proofs of
analytical inequalities remains a challenging and important task of
artificial intelligence and automated reasoning \cite{Bundy_1983},
\cite{Kaliszyk_Wiedijk_2007}.

\smallskip
Notice the logical-hardness general problem under consideration.
According to Wang \cite{Wang_1974}, for every function $G$ defined
by arithmetic operations and a composition over polynomials and sine
functions of the form $\sin \pi x$, there is a real number $r$ such
that the problem $G(r)=0$ is undecidable (see \cite{Poonen_2014}).
In 2003, M. Laczkovich \cite{Laczkovich_2003} proved that this
result can be derived if the function $G$ is defined in terms of the
functions $x,\sin x$ and $\sin (x\sin x^{n})$, $n=1,2,\ldots $
(without involving $\pi $). On the other hand, several algorithms
\cite{Tarski_1951}, \cite
{Muresan_2008}
and
\cite{Narkawicz_Munoz_Dutle_2015}
have been developed to determine the sign and the real zeroes of a given
polynomial, so that such problems can be considered as decidable
(see also \cite{Cutland_1980},
\cite{Poonen_2014}).

\smallskip
Let us denote by
\[
\displaystyle T_{n}^{\phi, a}(x) = \sum_{k=0}^{n} \displaystyle\frac{\phi^{(k)}(a)}{k!}%
(x-a)^{k}
\]
\noindent the Taylor polynomial of $n$-th degree associated to the function $%
\phi$ at a point $a$. Here, $\overline{T}_{n}^{\phi ,a}(x)$ and
$\underline{T}_{n}^{\phi ,a}(x)$
represent the Taylor polynomial of $n$-th degree associated to the function $%
\phi $ at a point $a$, in case \mbox{$\displaystyle T_{n}^{\phi,
a}(x) \geq \phi(x)$}, respective \mbox{$\displaystyle T_{n}^{\phi,
a}(x) \leq \phi(x)$}, for every $x\in \left( a,b\right) .$ We will
call $\overline{T}_{n}^{\phi ,a}(x)$ and $\underline{T}_{n}^{\phi
,a}(x) $ an upward, respective a downward approximation of $\phi ,$
on $\left( a,b\right) .$

\smallskip
We present a new algorithm for approximating a given MTP function
$f(x)$ by a polynomial function $P(x)$ such that
\begin{equation}
f(x)>P(x),  \label{N-A-goal}
\end{equation}%
using the upward and downward Taylor approximations $\underline{T}%
_{\,n}^{\,\sin ,0\!}{(x)}$, $\overline{T}_{n}^{\,\sin ,0\!}{(x)}$, $%
\underline{T}_{\,n}^{\,\cos ,0\!}{(x)}$, $\overline{T}_{n}^{\,\cos
,0\!}{(x)} $.

\section{The natural approach method and the associated algorithm}


The following two lemmas \cite{Malesevic_Makragic_2015}  related to
the Taylor polynomials associated to sine and cosine functions will
be of great help in our study.
\begin{lemma}
\label{S}
Let $T_{n}(x)\!=\!\!\mathop{\mbox{$\displaystyle\sum$}}\limits%
_{i=0}^{(n-1)/2}\displaystyle\frac{(-1)^{i}x^{2i+1}}{(2i+1)!}$.\newline
{$(i)$} If $n=4s+1$, with $s\in \NN_{0},$ then:%
\begin{equation}
\!\!T_{n}(x)\geq T_{n+4}(x)\geq \sin x,\,\,\,\, \mbox{for
every}\,\,\,\,0\leq x\leq \sqrt{(n+3)(n+4)};
\end{equation}%
and%
\begin{equation}
\!\!T_{n}(x)\leq T_{n+4}(x)\leq \sin x,\,\,\,\, \mbox{for every}\,\,-\sqrt{%
(n+3)(n+4)}\leq x\leq 0.
\end{equation}%
{$(ii)$} If $n=4s+3$, with $s\in \NN_{0},$ then:%
\begin{equation}
\!\!T_{n}(x)\leq T_{n+4}(x)\leq \sin x, \,\,\,\, \mbox{for
every}\,\,\,0\leq x\leq \sqrt{(n+3)(n+4)};
\end{equation}%
and%
\begin{equation}
\!\!T_{n}(x)\geq T_{n+4}(x)\geq \sin x, \,\,\,\,  \mbox{for every}\,\,-\sqrt{%
(n+3)(n+4)}\leq x\leq 0.
\end{equation}
\end{lemma}

\begin{lemma}
\label{C}
Let $T_{n}(x)=\!\mathop{\mbox{$\displaystyle\sum$}}\limits_{i=0}^{n/2}%
\displaystyle\frac{(-1)^{i}x^{2i}}{(2i)!}$.\newline
{$(i)$} If $n=4k$, with $k\in \mathbb{N}_{0}$, then:%
\begin{equation}
T_{n}(x)\geq T_{n+4}(x)\geq \cos x,\,\,\mbox{for every}\,\,-\!\sqrt{%
(n\!+\!3)(n\!+\!4)}\!\leq \!x\!\leq \!\sqrt{(n\!+\!3)(n\!+\!4)}.\!
\end{equation}%
{$(ii)$} If $n=4k+2$, with $k\in \mathbb{N}_{0}$, then:%
\begin{equation}
T_{n}(x)\leq T_{n+4}(x)\leq \cos x,\,\,\mbox{for every}\,\,-\!\sqrt{%
(n\!+\!3)(n\!+\!4)}\!\leq \!x\!\leq \!\sqrt{(n\!+\!3)(n\!+\!4)}.\!
\end{equation}
\end{lemma}

According to Lemmas~\ref{S}-\ref{C}, the upper bounds of the
approximation
intervals of the functions $\sin x$ and $\cos x$ are $\varepsilon _{1}=\sqrt{%
(n_{1}+3)(n_{1}+4)}$ and $\varepsilon
_{2}=\sqrt{(n_{2}+3)(n_{2}+4)}$, respectively. As $\varepsilon
_{1}>\mbox{\small$\displaystyle\frac{\pi}{2}$}$ and $\varepsilon
_{2}>\mbox{\small$\displaystyle\frac{\pi}{2}$}$, the
results of these lemmas are valid in particular, in the entire interval $\left( 0,%
\mbox{\small$\displaystyle\frac{\pi}{2}$}\right) $.
\begin{lemma} $\,$  \\
1) Let $n\in \NN $ and $x\in \left( 0,\mbox{\small
$\displaystyle\frac{\pi}{2}$}\right) $. Then:
\[
T_{n}^{\,\sin ,0}(x)\,\geq \,0.
\]%
2) Let $s\in \NN_{0}$, $\,p\in \NN $ and $x\in \left( 0,\mbox{\small
$\displaystyle\frac{\pi}{2}$}\right) $. Then:
\[
\left( \underline{T}_{4s+3}^{\,\sin ,0}(x)\right) ^{p}\leq \sin ^{p}{\!x}%
\leq \left( \overline{T}_{4s+1}^{\,\sin ,0}(x)\right) ^{p}.
\]
\end{lemma}
\begin{lemma}  $\,$
\label{Lemma_4} \\
Let $k\in
\mathbb{N}
_{0}$, $p\in
\mathbb{N}
$ and $x\in \left( 0,\mbox{\small
$\displaystyle\frac{\pi}{2}$}\right) $. Then:
\[
\cos ^{\,p}{\!x\leq }\left( \overline{T}_{4k}^{\,\cos ,0}(x)\right)
^{p}.
\]
\end{lemma}

In contrast to the function $\sin x$ and its downward Taylor
approximations, in the interval $\left( 0,\mbox{\small
$\displaystyle\frac{\pi}{2}$}\right) $ the function $\cos x$ and the
downward Taylor approximations $\underline{T}_{\,4k+2}^{\,\cos ,0}(x)=%
\mathop{\mbox{$\displaystyle\sum$}}_{i=0}^{2k+1}{%
\mbox{\small $\displaystyle\frac{(-1)^ix^{2i}}{(2i)!}$}}$, $\,k \in
\NN_{0}${,} require
special attention as there is no downward Taylor approximation $\displaystyle%
\underline{T}_{\,4k+2}^{\;\cos ,0}(x)$, such that $\cos
^{2}{\!x}\,\geq \,\left( \underline{T}_{\,4k+2}^{\;\cos
,0}(x)\right) ^{2},$ for every $x\in \left( 0,\mbox{\small
$\displaystyle\frac{\pi}{2}$}\right) $.


\smallskip
We present the following results related to the problem with downward
Taylor approximations of the cosine function.

\begin{proposition}
\label{cos-nanizni-polinomi}$\,$ \\
1) For every $k\in \mathbb{N}_{0}\,$, the downward Taylor approximation $%
\displaystyle\underline{T}_{\,4k+2}^{\;\cos ,0}(x)$ is a strictly
decreasing function on $\left( 0,\mbox{\small
$\displaystyle\frac{\pi}{2}$}\right) $.\newline 2) For every $k\in
\mathbb{N}_{0}\,$, there exists an unique $c_{k}\in
\left( 0,\mbox{\small $\displaystyle\frac{\pi}{2}$}\right) $ such that %
\mbox{$\,\displaystyle\underline{T}_{\,4k+2}^{\;\cos
,0}(c_{k})=0$.}\newline 3) The sequence $\displaystyle\left(
c_{k}\right) _{k\in
\mathbb{N}
_{0}},$ with $c_{0}=\sqrt{2}$, is strictly increasing and %
\mbox{$\displaystyle\,\lim_{k\rightarrow +\infty }{\!c_{k}}=\frac{\pi }{2}$.}%

\break

\noindent
4) For every $k\in  \NN_{0}$, there exists $d_{k}\in \left( c_{k},%
\mbox{\small
$\displaystyle\frac{\pi}{2}$}\right) $ such that $\cos {d_{k}}%
\,=\,\left\vert \,\underline{T}_{\,4k+2}^{\;\cos {},0}(d_{k})\,\right\vert $.%
\newline
5) The sequence $\displaystyle\left( d_{k}\right) _{k\in
\mathbb{N}
_{0}}$ is strictly increasing and $\displaystyle\,\lim_{k\rightarrow
+\infty }{\!d_{k}}=\frac{\pi }{2}$.
\end{proposition}

\noindent \textbf{Proof.} \textit{1)} The function $\displaystyle\underline{T%
}_{\,4k+2}^{\;\cos ,0}(x)$ is strictly decreasing on $\left( 0,%
\mbox{\small $\displaystyle\frac{\pi}{2}$}\right) ,$ since,
according to Lemma 1, $\,\, \displaystyle\left(
\,\underline{T}_{\,4k+2}^{\;\cos ,0}(x)\,\right) ^{\prime
}\,=\,-\overline{T}_{4k+1}^{\,\sin ,0}\,(x)\,\leq \,0. $

\medskip \noindent \textit{2)} The existence of $c_{k}$ follows from the
fact that $\displaystyle\underline{T}_{\,4k+2}^{\;\cos
,0}(0)\!=\!1\,>\,0$ and
$\displaystyle\,\underline{T}_{\,4k+2}^{\;\cos ,0}\left(
\mbox{\small $\displaystyle\frac{\pi}{2}$}\right) \,<\,\cos {\left(
\mbox{\small $\displaystyle\frac{\pi}{2}$}\right) }\!=\!0$.

\medskip \noindent \textit{3)} The monotonicity of the sequence $%
\displaystyle{\big (}c_{k}{\big )}_{k\in
\mathbb{N}
_{0}}$ is a result of the monotonicity of $\displaystyle\,\underline{T}%
_{\,4k+2}^{\;\cos ,0}(x)$ and \mbox{Lemma \ref{C} (ii).} \newline
The convergence of the sequence $\displaystyle{\big (}T_{n}^{\;\cos ,0}(x){%
\big )}_{n\in
\mathbb{N}
}$ implies the convergence of the sequence $\displaystyle{\big (}c_{k}{\big )%
}_{k\in \NN_{0}}$ to \mbox{\small $\displaystyle \,\frac{\pi}{2}$}.

\medskip \noindent \textit{4)} The function $\left| \,\underline{T}^{\,\cos{}%
,0}_{4k+2}(x) \,\right|$ is decreasing on $\displaystyle  (0,
c_{k})$ and increasing on $\displaystyle  \left(c_{k}, \mbox{\small
$\displaystyle\frac{\pi}{2}$}\right)$. Based on \mbox{Lemma \ref{C}
(ii)},  it follows that there exists $d_{k} \in \left( c_{k},
\mbox{\small
$\displaystyle\frac{\pi}{2}$} \right)$ such that $\cos{d_{k}}\, =\,\left|\,%
\underline{T}^{\;\cos{},0}_{\,4k+2}(d_{k})\,\right|$.

\medskip \noindent \textit{5)} This statement is a  consequence of the monotonicity of the sequence  $%
\displaystyle{\big (}c_{k}{\big )}_{k\in
\mathbb{N}
_{0}}$ and the increasing  monotonicity of  the function $\left\vert \,\underline{T}%
_{4k+2}^{\,\cos {},0}(x)\,\right\vert $  on $\displaystyle  \left(
c_{k},\mbox{\small $\displaystyle
\frac{\pi}{2}$}\right) $. \hfill $\blacksquare $ 

\begin{corollary} $\,$ \\
Let $\,k\in \NN_{0}$ and $\,p\in \NN $. Then:

\noindent 1) $\displaystyle\cos ^{2p}{\!x}>\left( \underline{T}%
_{4k+2}^{\,\cos ,0}(x)\right)^{2p}\!\!\!,\,$ for every
$\,\displaystyle x\in {\big (}0,d_{k}{\big )};$

\noindent 2) $\displaystyle\cos ^{2p}{\!x}<\left( \underline{T}%
_{4k+2}^{\,\cos ,0}(x)\right) ^{2p}\!\!\!,\,$ for every $x\in {\big (}d_{k},%
\mbox{\small $\displaystyle\frac{\pi}{2}$}{\big )}$.
\end{corollary}

\medskip
\noindent
Based on the above results, we have:

\begin{corollary} $\,$ \\
Let $\,k\in \NN_{0}$ and  $\,p\in \NN $. Then
$\underline{T}_{\,4k+2}^{\;\cos {},0}(x)$ is
not a downward approximation of the MTP function $\cos ^{2p}{\!x}\ $on $%
\displaystyle\left( d_{k},\mbox{\small
$\displaystyle\frac{\pi}{2}$}\right) $.
\end{corollary}

In order to ensure the correctness  of the algorithm (
\cite{Cutland_1980}, \cite{Algorithms}) we will  develop  next in
the sequel,  the following problem needs to be considered:

\medskip \smallskip

\noindent \textbf{Problem. }

\smallskip \noindent
\emph{For a given $\delta\!\in\!\left(0,\mbox{\small $\displaystyle\frac{\pi}{2}$}\right)$
and
$\,\mathcal{I}\!\subseteq\!\left(0,\mbox{\small $\displaystyle\frac{\pi}{2}$}\right)$,
find
$\displaystyle\widehat{k}\!\in\!\mathbb{N}_{0}$
such that for all
$\,k\!\in\!\mathbb{N}_{0}$, $k\!\geq\!\widehat{k}$ and $x\!\in\!\mathcal{I}\,:$}
\begin{equation}
\cos ^{2}{\!x}\geq \left( \underline{T}_{\,4k+2}^{\;\cos,0}(x)\right) ^{2}.
\end{equation}


\medskip
\noindent
\textbf{Remark.} \emph{If $\cos x$ appears in odd powers
only in the given MTP function $f(x)$, we take $\widehat{k}=0$}.

\break

\medskip
One of the method  to solve the problem of downward approximation of
the function $\cos ^{2p}{\!x},\,p\in
\mathbb{N}
$ is the {\bf method of multiple angles} developed in \cite%
{Malesevic_Makragic_2015}. All degrees of the functions $\sin x$ and
$\cos x$
are eliminated from the given MTP function $f(x)$%
, through conversion into multiple-angle expressions. This removes
all even degrees of the function $\cos x$, but then sine and cosine
functions appear in the form $\sin $\mbox{\boldmath $\kappa$}$x$ or $\cos $%
\mbox{\boldmath $\kappa$}$x$ where $\mbox{\boldmath $\kappa$}\,
x\!\in\!\left( 0,\mbox{\boldmath $\kappa$}\, \mbox{\small
$\displaystyle\frac{\pi}{2}$}\right)$ and $\mbox{\boldmath $\kappa$}
\!\in\! \mathbb{N}$. In this case, in order to use the results of
Lemmas~\ref{S}-\ref{C}, we are forced to choose large enough values
of $k \!\in\! \mathbb{N}_{0}$, such that $\sqrt{(k+3)(k+4)} \!>\!
\mbox{\boldmath $\kappa$}\,\mbox{\small
$\displaystyle\frac{\pi}{2}$}$. Note that higher value of $k$
implies a higher degree of the downward Taylor approximations and of
the polynomial $P(x)$ in (\ref{N-A-goal}) (for instance, see
\cite{Nenezic_Malesevic_Mortici_2015} and
\cite{Malesevic_Banjac_Jovovic_2015}).

\medskip
\noindent Several more ideas to solve the above problem are proposed
and considered below, under the names of Method A-D. In the
following, the numbers $c_{k}$ and $d_{k}$ are those defined in
Proposition \ref{cos-nanizni-polinomi}.

\bigskip \noindent \textbf{$\underline{\mbox{Method  A}}$}

\smallskip \noindent \emph{If $\displaystyle\delta \,<\,%
\mbox{\small $\displaystyle\frac{\pi}{2}$}$, find the smallest $k\in
\mathbb{N}
_{0}$ such that $\displaystyle d_{k}\in \left( \delta ,%
\mbox{\small $\displaystyle\frac{\pi}{2}$}\right) $. Then
$\widehat{k}=k$.}

\medskip \noindent Note that Method A assumes the solving of a
transcendental equation of the form $\displaystyle\cos
{x}=\underline{T}_{\,4k+2}^{\;\cos ,0}(x),$ that requires numerical
methods.

\bigskip \noindent \textbf{$\underline{\mbox{Method B}}$}

\smallskip \noindent \emph{If $\displaystyle\delta \,<\,%
\mbox{\small $\displaystyle\frac{\pi}{2}$}$, find the smallest $k\in
\mathbb{N}
_{0}$ such that $\displaystyle c_{k}\in \left( \delta ,%
\mbox{\small $\displaystyle\frac{\pi}{2}$}\right) $. Then
$\widehat{k}=k$.}

\bigskip \noindent \textbf{$\underline{\mbox{Method C}}$}

\smallskip \noindent \emph{If $\displaystyle\,\delta \,<\,%
\mbox{\small $\displaystyle\frac{\pi}{2}$}$, find the smallest $k\in
\mathbb{N}
_{0}$ such that $\displaystyle\,\underline{T}_{\,4k+2}^{\;\cos
,0}(\delta )\,\geq \,0$. Then $\widehat{k}=k$.}

\bigskip
\noindent
Note that Method B and Method C return the same output as for a given
$\delta $ and for every $k\in \NN_{0}$ the following
equivalence holds true:
\[
\left( \,\displaystyle c_{k}\in \left( \delta ,\displaystyle\frac{\pi }{2}%
\right) \,\,\wedge \,\,\underline{T}_{\,4k+2}^{\;\cos
,0}(c_{k})\,=\,0\,\right) \,\,\,\Longleftrightarrow \,\,\,\underline{T}%
_{\,4k+2}^{\;\cos ,0}(\delta )\,\geq \,0.
\]%
As Method B assumes the  determining the root $c_{k}$ of the
downward Taylor approximation
$\displaystyle\underline{T}_{\,4k+2}^{\;\cos ,0}(x)$ and Method C
assumes the checking the sign of the downward Taylor approximation
at point the $x=\delta $, it is notable that Method C presents a
faster and simpler procedure.

\bigskip \noindent \textbf{$\underline{\mbox{Method D}}$}

\smallskip \noindent \emph{Eliminate all even degrees of the function $\cos
x $ using the transformation%
\begin{equation}\label{transformation}
\cos ^{2p}{\!x}\,=\,\left( 1-\sin ^{2}{\!x}\right) ^{p}\,=\,%
\mathop{\displaystyle
\sum}_{i=0}^{p}(-1)^{i}%
\mbox{\footnotesize $\left(\!\!
\begin{array}{c}
\mbox{\footnotesize $\!p\!$} \\[0.5 ex]
\mbox{\footnotesize $\!i\!$}
\end{array}
\!\!\right)$}\sin ^{2i}{\!x}. \end{equation}
 }

\vspace*{-2.5 mm}

\noindent \emph{Then $\widehat{k} = 0$.}

\break

\noindent
Note that Method D can be applied for any $0<\delta \leq \pi /2$.
Hence, if a MTP function $f(x)$ is considered in the whole interval
$\left( 0,\mbox{\small $\displaystyle\frac{\pi}{2}$}\right)$,
then Method D is applicable only (apart from the multiple-angle method).
However, Method D implies an increasement of the number of terms needed
to be estimated. Let us represent a given MTP function $f$ in the
following form:
\begin{equation}
f(x)=\sum\limits_{i=1}^{\mbox{\scriptsize \sf m}}\alpha
_{i}x^{p_{i}}\!\cos ^{2k_{i}}\!x\sin ^{r_{i}}\!x\,\,+\,\,f_{1}(x)
\label{parni-cos}
\end{equation}%
where there are no terms of the form $\cos ^{2J}{\!x},\,\,J\!\in \!%
\mathbb{N}
,$ in $f_{1}(x)$. The elimination of all terms of the form $\cos ^{2k_{i}}{%
\!x}$ from (\ref{parni-cos}) using  the transformation
(\ref{transformation}),
will increase the number of addends in (\ref{parni-cos}), in the
general case with
\[
k_{1}+k_{2}+\ldots +k_{\mbox{\scriptsize \sf m}};
\]
consequently,  it will increase the number of terms of the form
$\sin ^{j}{\!x}$, $j\in
\mathbb{N}
$,
in (\ref{parni-cos}) needed to be estimated.


\subsection{An algorithm based on the natural approach method}


Let $f$ be a MTP function and $\mathcal{I} \subseteq \left( 0,\pi/2\right)$.
We concentrate to find a polynomial $\mathcal{TP}^{f}(x)$ such that
for every $x\in \mathcal{I},$
\[
\displaystyle f(x)>\mathcal{TP}^{f}(x).
\]%
\noindent In this case, the associated MTP inequality $\,f(x)>0$ can be
proved if we show that  for every $x\in \mathcal{I},$%
\[
\mathcal{TP}^{f}(x)>0,
\]%
\noindent which is a decidable problem according to Tarski
\cite{Tarski_1951}, \cite{Poonen_2014}.

\medskip
The following algorithm describes the method for finding such
a polynomial $\mathcal{TP}^{f}(x)$.

\renewcommand*{\thefootnote}{(\arabic{footnote})}

\medskip

$\underline{\mathbf{ALGORITHM ~  \emph{Natural Approach}}}$

\vspace{1mm} \textbf{INPUT:} ~~function $f$, ~~$\delta \in \left( 0,\frac{%
\pi }{2}\right) .$

\textbf{OUTPUT:} ~~polynomial $\mathcal{TP}^{f}(x).$


{\footnotesize

\hspace*{3mm} \textbf{1.} $/\ast $  Solve a problem
involving downward approximations depending on $\cos ^{2p}{\!x},$
$\ast / $

\hspace*{10mm} $/\ast $  {\footnotesize i.e. determining $\widehat{k}\in \NN%
_{0}$, such that for all $k\in \NN_{0},$ $k\geq \widehat{k},$ it holds: }$%
\ast /$

\hspace*{9mm} $/\ast $  {\footnotesize $\,\,\,\displaystyle\cos
^{2}x\,\geq \,\left( \,\underline{T}_{\,4k+2}^{\;\cos
,0}(x)\,\right) ^{2},$ \ for\thinspace every \thinspace $x\in \left(
0,\delta \right] .\,\,$}~ $\ast / $

\hspace*{7mm} \textbf{If} ~~$\delta > \sqrt{2}$ and there are even degrees
of the function $\cos{x}$ \textbf{then}

\hspace*{25mm} \textbf{If} ~~$\displaystyle\delta <\mbox{\small
$\displaystyle\frac{\pi}{2}$}$~~\textbf{then} ~~ use $\underline{\mbox{Method C}}$ or $%
\underline{\mbox{Method D}}$

\hspace*{25mm} \textbf{else} ~~use $\underline{\mbox{Method D}}$

\hspace*{8mm} \textbf{else} ~ $\widehat{k} := 0$

\vspace{2mm} \hspace*{3mm} \textbf{2.} $/\ast $  {\footnotesize In
the procedure \emph{Estimation}} {\footnotesize (described below),
for a given MTP function } $\ast /$

\hspace*{9mm} $/\ast $  {\footnotesize $f(x),$ each addend
$a_{i}(x)$ in the function $f(x)$ is estimated.}~ $\ast /$

\textbf{PROCEDURE Estimation}\textbf{\boldmath $\,\left( \, f(x) \, \right)$}

\smallskip \textbf{END} /$\ast$ Algorithm $\ast$/

}

\break

\bigskip $\underline{\mathbf{PROCEDURE ~Estimation }}$

\vspace{1mm} \textbf{INPUT:} ~~the function
$f(x)=\sum_{i=1}^{n}a_{i}(x),\,$ where $a_{i}(x)\,=\,\alpha
_{i}x^{p_{i}}\!\cos ^{q_{i}}{\!x}\sin ^{r_{i}}{\!x}$.

\textbf{OUTPUT:}
\begin{minipage}[t]{97mm} the
polynomial  ${\cal TP}^{f}(x)$ and array $\left(\,\left(s_{i}, k_{i}
\right)\,\right),\,i=1, ..., n$ where $s_{i}$ and $k_{i}$ represent
the number that determines the degree of the Taylor approximation of
the function $\sin{x}$, respective $\cos{x}$ in the addend
$a_{i}(x)$.
\end{minipage}

{\footnotesize

\smallskip Estimate each addend $a_{i}(x)$ with $q_{i}^{2} + r_{i}^{2} \neq 0$,  as follows:

\hspace*{12mm} \textbf{I} ~ If $\alpha_{i}>0$, then:

\hspace*{18mm} /$\ast$ First select the degrees of the downward
approximations $\ast$/

\hspace*{18mm} Select $s_{i}\geq 0$ and $k_{i} \geq \widehat{k}$.

\hspace*{18mm} Estimate: $a_{i}(x)\geq \alpha _{i}x^{p_{i}}\left( \underline{%
T}_{\,4s_{i}+3}^{\;\sin ,0}(x)\right) ^{q_{i}}\left( \underline{T}%
_{\,4k_{i}+2}^{\;\cos ,0}(x)\right) ^{r_{i}};$

\hspace*{12mm} \textbf{II} ~ If $\alpha _{i}<0$, (i.e. $\alpha
_{i}=-\beta _{i},$ with  $\beta _{i}>0$)

\hspace*{18mm} /$\ast$ First select the degrees of the downward
approximations $\ast$/

\hspace*{18mm} Select $s_{i}\geq 0$ and $k_{i}\geq 0$.

\hspace*{18mm} Estimate:

\hspace*{18mm} $~ a_{i}(x) =-\beta_{i}x^{p_{i}}(\sin x)^{q_{i}}(\cos
x)^{r_{i}}\geq - \beta_{i}x^{p_{i}}\left(\overline{T}^{\,%
\sin,0}_{4s_{i}+1}(x)\right)^{q_{i}}\left(\overline{T}^{\,%
\cos,0}_{4k_{i}}(x)\right)^{r_{i}} $;

\vspace*{-1.0 mm}

\hspace*{10mm} $/\ast$ {\small Estimation of each addend $a_{i}(x)$
in function $f(x)$ yields a polynomial }$\ast/$

\hspace*{10mm}  $/\ast$ {\small  of the form:   }  $\ast/$

{\small \hspace*{9mm} $/\ast$ $\;P(x)=\mathop{\mbox{$\displaystyle \sum$}}%
_{i=1}^{n}
\alpha_{i}x^{p_{i}}\left(T_{\,n_{i}}^{\;\sin,0}(x)\right)^{q_{i}}\left(T_{%
\,m_{i}}^{\;\cos,0}(x)\right)^{r_{i}}$, where \mbox{\small $T \in \lbrace
\underline{T}, \overline{ T}\rbrace$}}. ~ $\ast/$

\hspace*{2mm} \textbf{Return:} the polynomial $\,\mathcal{TP}^{f}(x)
\,\,\,$ and \,\,\, array $\left( \,\left( s_{i}, k_{i} \right)
\,\right)$.

\smallskip  \textbf{END} /$\ast$ Procedure $\ast$/

}

\smallskip
\noindent \textbf{Comment} on step \textbf{II} of the Procedure
Estimation: in the general case, the addend $\, a_{i}(x)=\,\,-\beta
_{i}x^{p_{i}}(\sin x)^{q_{i}}(\cos x)^{r_{i}} $ can be estimated in
one of the following three ways:


$
\!\!\!\!\;\,(i)~~ a_{i}(x) =-\beta_{i}x^{p_{i}}(\sin x)^{q_{i}}(\cos
x)^{r_{i}}\geq \beta_{i}x^{p_{i}}\left(\underline{T}^{\,%
\sin,0}_{4s_{i}+3}(x)\right)^{\!q_{i}}\!\!\left(-\overline{T}^{\,%
\cos,0}_{4k_{i}}(x)\right)^{\!r_{i}}\!\!\!\!,$


$
\!\!\!\!\,(ii)~~ a_{i}(x) =-\beta_{i}x^{p_{i}}(\sin x)^{q_{i}}(\cos
x)^{r_{i}}\geq\beta_{i}x^{p_{i}}\left(-\overline{T}^{\,\sin,0}_{4s_{i}+1}(x)%
\right)^{\!q_{i}}\!\!\left(\underline{T}^{\,\cos,0}_{4k_{i}+2}(x)\right)^{\!r_{i}}\!\!\!\!,
$


$
\!\!\!\!(iii)~~ a_{i}(x) =-\beta_{i}x^{p_{i}}(\sin x)^{q_{i}}(\cos
x)^{r_{i}}\geq - \beta_{i}x^{p_{i}}\left(\overline{T}^{\,%
\sin,0}_{4s_{i}+1}(x)\right)^{\!q_{i}}\!\!\left(\overline{T}^{\,%
\cos,0}_{4k_{i}}(x)\right)^{\!r_{i}}\!\!\!\!.$

\medskip
\noindent Note that for fixed $s_{i}, k_{i}, q_{i} $ and $r_{i}$,
the  method $(iii)$ generates  polynomials of the smallest degree.

\noindent We present the  following characteristic (\cite{Knuth},
\cite{Algorithms}) for  the \emph{Natural Approach} algorithm.

\begin{theorem}
The Natural Approach algorithm  is correct.
\end{theorem}

\noindent \textbf{Proof.}  Every step in the algorithm is based  on
the results  obtained from  Lemmas \ref{S}-\ref{Lemma_4} and
Proposition \ref{cos-nanizni-polinomi}. Hence, for every input
instance (i.e.$\;$for any  MTP function $f(x)$ over a given interval
${\cal I}\subseteq \left(0, \pi/2\right)$), the algorithm halts with
the correct output (i.e. the algorithm returns
the corresponding polynomial).\hfill $\blacksquare $ 


\section{Some applications of the algorithm}


\smallskip

We present  an application of the \emph{Natural Approach} algorithm
in the proof (Application 1  -  Theorem \ref{Tanjina-teorema}) of
certain new rational (Pad$\acute{\mbox{e}}$) approximations of the
function $\cos ^{2}{\!x}$, as well as in the improvement  of  a
class of inequalities (\ref{Z-H-Jang-1}) by Z. H. Yang (Application
2, Theorem \ref{Brankova-teorema}).

\break

\medskip \noindent \textbf{Application 1}

\medskip
\noindent Bercu \cite{Bercu_2016} used the  Pad$\acute{\mbox{e}}$
approximations  to prove certain inequalities for trigonometric
functions. Let us denote by $\left( f(x)\right) _{[m/n]}$ the
Pad$\acute{\mbox{e}}$ approximant $\left[ m/n\right] $ of the
function $f(x)$.

\medskip \noindent In this example we introduce a constraint of the function
$\cos^{2}{\!x}$ by the following Pad$\acute{\mbox{e}}$
approximations:
\[
\displaystyle
\left(\cos^{2}{\!x}\right)_{[6/4]} = \displaystyle\frac{-59\,x^{6} + 962\,x^{4} -
3675\,x^{2} + 4095}{17\,x^{4} + 420\,x^{2}+ 4095}
\]
and
\[
\left(\cos^{2}{\!x}\right)_{[4/4]} = \displaystyle\frac{163\,x^{4} - 780\,x^{2} + 945}{%
13\,x^{4} + 165\,x^{2} + 945}.
\]

\begin{theorem}
\label{Tanjina-teorema} The following inequalities hold true, for every $%
\displaystyle x \!\in\! \left( 0,\mbox{\small $\displaystyle\frac{\pi}{2}$}\right)\!:$
\begin{equation}
\displaystyle \left(\cos^{2}{\!x}\right)_{[6/4]} %
\,<\,\cos ^{2}{\!x}\,<\,\displaystyle
\left(\cos^{2}{\!x}\right)_{[4/4]} \label{Tanjin-primer}
\end{equation}
\end{theorem}

\noindent \textbf{Proof.} We  first prove \textbf{the left-hand
side} inequality (11). Using a computer software for symbolic
computations, we can
conclude that the  function $G_{1}(x) \,=\,\left( \cos ^{2}{\!x}%
\right) _{[6/4]}\,\ $has exactly one zero $\delta =1.551413...$ in the
interval $\left( 0,\,\mbox{\small $\displaystyle\frac{\pi}{2}$}\right) .$ As $%
G_{1}(0)=1>0$ and $G_{1}\left( \mbox{\small $\displaystyle\frac{\pi}{2}$}\right)
=-0.000431...\,<0$, we deduce that%
\begin{equation}
G_{1}(x)\geq 0\;\;\;\mbox{for every}\;\;\;x\in (0,\delta ]
\label{levo-od-delta}
\end{equation}%
and
\begin{equation}\label{desno-od-delta}
G_{1}(x)<0\;\;\;\mbox{for every}\;\;\;x\in \left(\delta ,\mbox{\small
$\displaystyle\frac{\pi}{2}$}\right).
\end{equation}%
Moreover, $G_{1}(x)<\cos ^{2}{\!x}$, for every $x\in \left( \delta ,\,%
\mbox{\small
$\displaystyle\frac{\pi}{2}$}\right) $. We prove now that%
\begin{equation}
\displaystyle G_{1}(x)<\cos ^{2}{\!x},\,\,\,x\in \left( 0,\,\delta \right] .
\label{levo}
\end{equation}%
We search a downward Taylor polynomial $\underline{T}_{\,4k+2}^{\;\cos ,0}(x)
$, such that for every $x\in \left( 0,\,\delta \right] $
\begin{equation}
\displaystyle G_{1}(x)<\left( \underline{T}_{\,4k+2}^{\;\cos ,0}(x)\right)
^{2}<\cos ^{2}{\!x}.  \label{interpolacija}
\end{equation}

\medskip \noindent We apply the \emph{Natural Approach} algorithm to the
function $f(x)=\cos ^{2}{\!x},\,$ $x\in (0,\delta ]$, to determine the
downward Taylor polynomial $\underline{T}_{\,4k+2}^{\;\cos ,0}(x)$, such
that
\[
\displaystyle
\left(
\underline{T}_{4k+2}^{\,\cos ,0}(x)
\right)^{2}
\!<
\cos ^{2}{\!x,}
\;\;x\in (0,\delta].
\]%
We can use Method C, or Method D from the \emph{Natural Approach} algorithm, since $%
\delta <\mbox{\small $\displaystyle\frac{\pi}{2}$}.$ In this proof,
we choose Method C.

\setcounter{footnote}{0}

\smallskip
\noindent The smallest $k$ for which $\underline{T}_{\,4k+2}^{\;\cos
,0}(\delta )>0$ is $k=1$. Therefore $\widehat{k}=1$. In the
\emph{Estimation } procedure only step I can be applied to the
(single)
addend $\cos ^{2}{\!x}$. In this step, $s_{1}\geq 0$ and $k_{1}\geq \widehat{%
k}=1$ should be selected. Let us select $s_{1}=0$ and $k_{1}=2.$\footnote{%
For the selection $s_{1}=0$ and $k_{1}=1$, the output of the \emph{Natural
Approach} algorithm is the polynomial:
\[
\mathcal{TP}(x)=\underline{T}_{6}^{\,\cos ,0}(x)=  1-\displaystyle\frac{x^{2}}{2!}+\displaystyle\frac{%
x^{4}}{4!}-\displaystyle\frac{x^{6}}{6!}
\]%
such that $\mathcal{TP}(x)\!\!%
\begin{array}{c}
\mbox{\scriptsize $<$} \\[-1.25ex]
\mbox{\scriptsize $>$}%
\end{array}%
\!\!G_{1}(x)$ holds for some $x\in (0,\delta ]$.} As a result of
this selection, the output of the \emph{Natural Approach} algorithm
is the polynomial:
\[
\mathcal{TP}(x)=\left( \underline{T}_{10}^{\,\cos ,0}(x)\right) ^{2}=\left(1-%
\displaystyle\frac{x^{2}}{2!}+\displaystyle\frac{x^{4}}{4!}-\displaystyle\frac{x^{6}}{6!}+\displaystyle\frac{x^{8}}{8!}-%
\displaystyle\frac{x^{10}}{10!}\right)^{2}.
\]%
We prove that
\begin{equation}
\left( \underline{T}_{10}^{\,\cos ,0}(x)\right) ^{2}-G_{1}(x)>0,\;\;
x\in \left( 0,\,\delta \right] .  \label{levo-dokaz}
\end{equation}%
This is true, since
\[
\left( \underline{T}_{10}^{\,\cos ,0}(x)\right) ^{2}-G_{1}(x)=\displaystyle\frac{x^{12}}{%
13168189440000\,(17\,x^{4}+420\,x^{2}+4095)}\,Q(x),
\]%
where
\[
\begin{array}{rcl}
Q(x)\! & \!=\! & \!17\,x^{12}+15\,x^{8}(15837-176\,x^{2})+8100%
\,x^{4}(64519-1687\,x^{2}) \\[2ex]
\! & \!\! & \!+3200\,(50205015-4035906\,x^{2})\,>\,0.%
\end{array}%
\]%
Finally, we have $\,G_{1}(x)<\cos ^{2}{\!x}$ for every $x\in \left(
0,\,\delta \right] $. According to (\ref{desno-od-delta}), we have
\[
G_{1}(x)<\cos ^{2}{\!x}, \;\;\;\mbox{for every}\;\;\;x\in \left( 0,\,%
\mbox{\small $\displaystyle\frac{\pi}{2}$}\right) .
\]%
Now we prove \textbf{the right-hand side} inequality
(\ref{Tanjin-primer}). For \ $G_{2}(x)=\left( \cos ^{2}{\!x}\right)
_{[4/4]}$ we prove the following inequalities, for every $x\in
\left( 0,\mbox{\small $\displaystyle\frac{\pi}{2}$}\right) $:
\begin{equation}
\displaystyle\cos ^{2}{\!x}<\left( \overline{T}_{8}^{\,\cos ,0}(x)\right)
^{2}<G_{2}(x).  \label{interpolacija-2}
\end{equation}%
Based on Proposition \ref{cos-nanizni-polinomi}, it is enough to prove that for
every $x\in \left(0,\mbox{\small $\displaystyle\frac{\pi}{2}$}\right)$
\begin{equation}
\left( \overline{T}_{8}^{\,\cos {x},0}(x)\right) ^{2}<G_{2}(x).
\label{dovoljno-za-desno}
\end{equation}%
This is true, as
\[
G_{2}(x)-\left( \overline{T}_{8}^{\,\cos {x},0}(x)\right) ^{2}\,= \displaystyle\frac{%
x^{10}\,}{1625702400\,(13\,x^{4}+165\,x^{2}+945)}\,R(x),
\]%
where%
\[
R(x)=x^{8}(1291-13x^{2})+x^{4}(2004240-66913x^{2})+480(632604-74625x^{2})\,>
\,0.\]%
Since $\cos ^{2}{\!x}\leq \left( \overline{T}_{4k}^{\,\cos ,0}(x)\right) ^{2}
$, for every $k\in \NN_{0}$ and all $x\in (0,\mbox{\small $\displaystyle\frac{\pi}{2}$})$%
, we have
\[
\cos ^{2}{\!x}<G_{2}(x),~~ \mbox{for every} \, x\in \left(
0,\,\mbox{\small $\displaystyle\frac{\pi}{2}$}\right) .
\]%
\hfill $\blacksquare $

\medskip \noindent \textbf{Note:} Using Pad$\acute{\mbox{e}}$
approximations, Bercu \cite{Bercu_2016}, \cite{Bercu_2017}
recently refined certain trigonometric inequalities over
various intervals
$\mathcal{I}=(0,\delta )
\subseteq
(0,  \mbox{\small $\displaystyle\frac{\pi}{2}$})$.
All such inequalities can be proved in a similar way and using
the \emph{natural approach} algorithm, as in the proof of
Theorem~\ref{Tanjina-teorema}.

\medskip
\medskip \noindent \textbf{Application 2.}

\medskip \noindent Z.-H. Jang \cite{Z.-H._Yang_2014} proved the following
inequalities, for every $x\in \left( 0,\pi \right) :$%
\begin{equation}
\cos ^{2}{\!\displaystyle\frac{x}{2}}\leq \displaystyle\frac{\sin {x}}{x}\leq \cos ^{3}{\!\displaystyle\frac{x}{3%
}}\leq \displaystyle\frac{2+\cos {x}}{3}.  \label{Z-H-Jang-1}
\end{equation}%
Previously, Kl$\acute{\mbox{\rm e}}$n, Visuri, and Vuorinen \cite%
{Klen_Visuri_Vuorinen_2010} proved the above inequality on ${\big (}0,\sqrt{%
27/5}\,{\big )}$ only.

\medskip In this example we propose the following improvement of (\ref%
{Z-H-Jang-1}):

\begin{theorem}
\label{Brankova-teorema} The following inequalities hold true, for every $x
\in \left(0, \pi\right)$ and $a \in \displaystyle
\left(1, \mbox{\small $\displaystyle\frac{3}{2}$}\right)$:
\begin{equation}  \label{Z-H-Jang}
\cos^{2}{\!\displaystyle\frac{x}{2}} \leq
\left(\displaystyle\frac{\sin{x}}{x}\right)^{a} \leq
\displaystyle\frac{\sin{x}}{x}.
\end{equation}
\end{theorem}

\noindent \textbf{Proof.} As $\displaystyle a>1$ and $\displaystyle0<\displaystyle\frac{%
\sin {x}}{x}<1$, we have:
\[
\left( \displaystyle\frac{\sin {x}}{x}\right) ^{\!a}\!<\displaystyle\frac{\sin {x}}{x}.
\]%
We prove now the following inequality:
\begin{equation}
\cos ^{2}{\!\displaystyle\frac{x}{2}}<\left( \displaystyle\frac{\sin {x}}{x}\right) ^{\!a}\!
\label{primer-1}
\end{equation}%
for every $x\in \left( 0,\pi \right) $ and $a\in \displaystyle\left( 1,%
\mbox{\small $\displaystyle\frac{3}{2}$}\right) $. It suffices to show that the
following mixed logarithmic-trigonometric-polynomial function \cite%
{Malesevic_Lutovac_Banjac_2015}
\begin{equation}
F(x)=a\ln \left( \displaystyle\frac{\sin {x}}{x}\right) -2\ln \left( \cos {\!\displaystyle\frac{x}{2}%
}\right)
\end{equation}%
is positive, for every $x\in \left( 0,\pi \right) $ and $a\in \left( 1,%
\mbox{\small $\displaystyle\frac{3}{2}$}\right) $. Given that
\begin{equation}
\lim_{x\rightarrow 0}F(x)=0,  \label{Limes 1}
\end{equation}%
based on the ideas from \cite{Malesevic_Lutovac_Banjac_2015}, we
connect the function $\,F(x)$ to the analysis of its derivative:
\[
\displaystyle F^{\prime\!}(x)=\displaystyle\frac{1}{2}\,\displaystyle\frac{
f\!\left(\mbox{\small $\displaystyle\frac{x}{2}$}\right)}{x\sin {%
\mbox{\small $\displaystyle\frac{x}{2}$}}\cos {\mbox{\small
$\displaystyle\frac{x}{2}$}}},
\]

\noindent where
\begin{equation}
f(t)=4t(a-1)\cos ^{2}{\!t}-2a\sin {t}\cos {t}-2t(a-2).
\end{equation}%
Let us note that $F^{\prime\!}(x)$ is the quotient of two MTP functions.

The inequality $\displaystyle F^{\prime\!}(x)>0$ is equivalent to $%
\displaystyle f(t)>0$. The proof of the later inequality will be done using
the \emph{Natural Approach} algorithm for the function $f(t)$ on $\left( 0,%
\mbox{\small $\displaystyle\frac{\pi}{2}$}\right) $, with $a\!\in \!\left( 1,%
\mbox{\small $\displaystyle\frac{3}{2}$}\right) $. As before, we
search a polynomial $\mathcal{TP}(t)$ such that
\[
f(t)>\mathcal{TP}(t)>0.
\]%
In the step 1 of the \emph{Natural Approach} algorithm, we can use
Method D only, because $\delta =\mbox{\small
$\displaystyle\frac{\pi}{2}$}$. Then
\begin{equation}
\begin{array}{rcl}
f(t) & = & 4\,t\,(a-1)\,(1-\sin ^{2}{t})\,\,-\,\,2\,a\,\sin {t}\,\cos {t}%
\,\,-\,\,2\,t\,(a-2) \\[2.5ex]
& = & 4\,t\,(1-a)\,\sin ^{2}{t}\,\,-\,\,2\,a\,\sin {t}\,\cos {t}%
\,\,+\,\,2\,t\,a%
\end{array}
\label{D}
\end{equation}%
with $\widehat{k}\,=\,0$.  In the \emph{Estimation} procedure only%
\footnote{%
Because for every fixed $a\in \left( 1,\frac{3}{2}\right) $: $\alpha
_{1}=4(1-a)<0$ and $\alpha _{2}=-2a<0$.} the step II can be applied
to the first and second addends in (\ref{D}), where $\displaystyle
s_{i}\geq 0$ and $\displaystyle k_{i}\geq 0$, $i=1,2$ should be
selected. Let us, for example, select $\displaystyle
s_{1}=k_{1}=s_{2}=k_{2}=1$. As a result of this selection, the
\emph{Natural Approach} algorithm yields the polynomial
\[
\mathcal{TP}(t)
=
4t(1\!-\!a)
\!\left(
t\!-\!
\mbox{\small $\displaystyle\frac{1}{6}$}t^{3}
\!+\!
\mbox{\small $\displaystyle\frac{1}{120}$}t^{5}
\right)^{\!2}
\!\!-
2a\!
\left(
t
\!-\!
\mbox{\small $\displaystyle\frac{1}{6}$}t^{3}
\!+\!
\mbox{\small $\displaystyle\frac{1}{120}$}t^{5}
\right)
\!\left(
1\!-\!\mbox{\small $\displaystyle\frac{1}{2}$}t^{2}
\!+\!
\mbox{\small $\displaystyle\frac{1}{24}$}t^{4}
\right)
\!+\!
2\,t\,a\,
\]%
for which $f(t)>\mathcal{TP}(t)$,  for every $t\!\in \!\left( 0,%
\mbox{\small
$\displaystyle\frac{\pi}{2}$}\right) $ and  $a\!\in \!\left( 1,%
\mbox{\small
$\displaystyle\frac{3}{2}$}\right) $. The inequality $f(t)>0$ is reduced to a decidable
problem:
\begin{equation}
\label{problem-PIP}
\mathcal{TP}(t)>0,\;\;\;\mbox{for every }\;\;\;t\in \left( 0,%
\mbox{\small $\displaystyle\frac{\pi}{2}$}\right) \,\,\,\mbox{and}\,\,\,a\in \left( 0,%
\mbox{\small $\mbox{\small $\displaystyle\frac{3}{2}$}$}\right) .
\end{equation}%
The sign of the polynomial $\mathcal{TP}(t)$ can be determined in several
ways. For example, let us represent the polynomial $\mathcal{TP}(t)$ as
\begin{equation}
\label{PIP-1}
\mathcal{TP}(t)\,=\,p(t)a\,+\,q(t),
\end{equation}%
where
\[
p(t)\!=\!-\displaystyle\frac{t^{3}\left(
2t^{8}\!-\!75t^{6}\!+\!1120t^{4}\!-\!7680t^{2}\!+\!19200\right) }{7200}
\]
and
\[
q(t)\!=\!4\,t\,{\Big (}t\!-\!\mbox{\small $\displaystyle\frac{1}{6}$}t^{3}\!+\! %
\mbox{\small $\displaystyle\frac{1}{120}$}t^{5}{\Big )}^{\!2}\!.
\]%
For every fixed $\displaystyle t\in \left( 0,\mbox{\small $\displaystyle\frac{\pi}{2}$}%
\right) $, the function $\displaystyle\mathcal{TP}(t)=p(t)a+q(t)$ is
 linear, monotonically  decreasing with respect to
$a\!\in\!\left(1,\mbox{\small $\displaystyle\frac{3}{2}$}\right)$,
since  for every \mbox{$t\!\in\!\left( 0,\mbox{\small
$\displaystyle\frac{\pi}{2}$}\right)$},
\[
\displaystyle p(t)=-\displaystyle\frac{t^{3}}{7200}\,{\Big
(}2\,t^{8}+5\,t^{4}(224-15t^{2})+3840\,(5-2t^{2}){\Big )}\,\,<0.
\]%
Hence, for every fixed $\displaystyle t\!\in \!\left( 0,\mbox{\small
$\displaystyle\frac{\pi}{2}$}\right) $, the value of (\ref{PIP-1})
is greater than the value of the same expression for $a=\mbox{\small
$\displaystyle\frac{3}{2}$}$:
\[
p(t)\,\mbox{\small $\displaystyle\frac{3}{2}$}\,+\,q(t)=-\displaystyle\frac{t^{5}}{14400}%
\,(2\,t^{6}-65\,t^{4}+800\,t^{2}-3840)
\]%
But
\[
p(t)\,\displaystyle\frac{3}{2}\,+\,q(t)=\displaystyle\frac{t^{5}}{14400}\,{\Big (}%
t^{4}(65-2\,t^{2})\,+160(\,24-5t^{2}){\Big )}>0,%
\]%
so the inequality (\ref{problem-PIP}) is true and consequently, $%
F^{\prime\!}(x)>0$ on $(0,\pi ),$ for every $a\in \left( 1,\mbox{\small
$\displaystyle\frac{3}{2}$}\right) $. But
$\displaystyle\lim_{x\rightarrow 0}{F(x)}=0$, so $F(x)>0$ on $(0,\pi)$,
for every $a\in \left( 1,\mbox{\small $\displaystyle\frac{3}{2}$}\right)$.

\smallskip
$\,\hfill \blacksquare $

\medskip
\noindent
\textbf{Remark on Theorem~\ref{Brankova-teorema}.}

\medskip \noindent Let us consider possible refinements of the inequality (\ref%
{Z-H-Jang-1}) by a real analytical function  $\varphi
_{a}(x)\!=\!\left( \displaystyle\frac{\sin x}{x}\right) ^{\!a}\!\!,
$ for $x\!\in
\!\left( 0,\delta \right) $ and $a\!\in \!%
\mathbb{R}
$. The function $\varphi _{a}(x)$ is real analytical, as it is related to
the analytical function%
\begin{equation}
t(x)=a\ln \!\left( \displaystyle\frac{\sin x}{x}\right) =a\mathop{\mbox{\large
$\displaystyle\sum$}}\limits_{k=1}^{\infty }{\mbox{\small
$\displaystyle\frac{(-1)^k2^{2k-1}B_{2k}}{k(2k)!}$}\,x^{2k}}
\end{equation}%
($B_{i}$ are the Bernoulli numbers, see e.g.
\cite{Gradshteyn_Ryzhik_2014}). The following consideration of the
sign of the analytical function in the left
and right neighborhood of zero is based on Theorem 2.5 from \cite%
{Malesevic_Makragic_2015}. Let us consider the real analytical function
\begin{equation}
f_{1}(x)
=
\left( \displaystyle\frac{\sin x}{x}\right) ^{\!a}-\cos ^{2}{\!\displaystyle\frac{x}{2}}
=\left(
-\mbox{\small $\displaystyle\frac{a}{6}$}+\mbox{\small $\displaystyle\frac{1}{4}$}
\right)x^{2}
+
\left(
\mbox{\small $\displaystyle\frac{{a}^{2}}{72}$}
-
\mbox{\small $\displaystyle\frac{a}{180}$}
-
\mbox{\small $\displaystyle\frac{1}{48}$}
\right){x}^{4}
+
\ldots,
\end{equation}
$x \in (0,\pi)$. The restriction
\begin{equation}
f_{1}^{\prime \prime }(0)=-\mbox{\small $\displaystyle\frac{a}{3}$}+\mbox{\small
$\displaystyle\frac{1}{2}$}>0
\end{equation}%
i.e.
\begin{equation}
a\in \left( -\infty ,\mbox{\small $\displaystyle\frac{3}{2}$}\right)
\end{equation}%

\break

\noindent
is a necessary and sufficient condition for $f_{1}(x)>0$ to hold on an
interval ${\big (}0,\delta _{1}^{(a)}{\big )}$ {\big (}for some $\delta
_{1}^{(a)}\!>\!0${\big )}. Also, the restriction
\begin{equation}
a\in \left( \mbox{\small $\displaystyle\frac{3}{2}$},\infty \right)   \label{uslov-1}
\end{equation}%
is a necessary and sufficient condition for $f_{1}(x)<0$ to hold on
an interval ${\big (}0,\delta _{2}^{(a)}{\big )}$ {\big (}for some
$\delta _{2}^{(a)}\!>\!0${\big )}. The following equivalences hold
true for every $x\!\in \!(0,\pi )$:
\begin{equation}
a\in \left( 1,\infty \right) \quad \Longleftrightarrow \quad \left( \displaystyle\frac{%
\sin x}{x}\right) ^{\!a}<\displaystyle\frac{\sin x}{x},  \label{uslov-2}
\end{equation}%
\begin{equation}
a\!\in \!\left( -\infty ,1\right) \quad \Longleftrightarrow \quad \displaystyle\frac{%
\sin x}{x}<\left( \displaystyle\frac{\sin x}{x}\right) ^{\!a}.
\end{equation}%
The refinement in Theorem~\ref{Brankova-teorema} is given based on
the possible values of the parameter $a$ in (33) and (34). A similar
analysis shows us that only the following refinements of the
inequality (\ref{Z-H-Jang-1}) are possible:

\begin{corollary}
Let $a\in \left[ \mbox{\small $\displaystyle\frac{3}{2}$},\,+\infty
\right) $. There exists $\delta >0$ such that for every $x\in \left(
0,\delta \right) $, it
holds:%
\begin{equation}
\left( \displaystyle\frac{\sin {x}}{x}\right) ^{a}\leq \cos ^{2}{\!\displaystyle\frac{x}{2}.}
\end{equation}
\end{corollary}

\begin{corollary}
Let $a\in \displaystyle\left( -\infty ,1\right) $. There exists $\delta >0$
such that for every $x\in \left( 0,\delta \right) $, it holds:
\begin{equation}
\displaystyle\frac{2+\cos {x}}{3}\leq \left( \displaystyle\frac{\sin {x}}{x}\right) ^{a}.
\end{equation}
\end{corollary}


\section{Conclusions and Future Work}


The results of our analysis could be implemented by means of an
automated proof assistant \cite{Geuvers}, so our work is a
contribution to the library of automatic support tools
\cite{Miller_2011} for proving various analytic inequalities.

\smallskip
Our general algorithm associated to the natural approach method can
be successfully applied to prove a wide category of classical MTP
inequalities.  For example, the \emph{Natural Approach} algorithm
has recently been used to prove some several open problems that
involve MTP inequalities (see e.g. \cite{Malesevic_Makragic_2015} -
\cite{Malesevic_Banjac_Jovovic_2015}).

\smallskip
It is our contention that  the \emph{Natural Approach} algorithm
can be used to introduce and solve other new similar results. Chen
\cite{Chen_2012} used a similar method to prove the following
inequalities,
for every $x\in \left( 0,1\right) $:%
\[
2+\frac{17}{45}x^{3}\arctan x<\left( \frac{\arcsin x}{x}\right)^{\!2}
\!+
\frac{\arctan x}{x}
\]%
and%
\[
2+\frac{7}{20}x^{3}\arctan x<2\left( \frac{\arcsin x}{x}\right)
\!+
\frac{\arctan x}{x};
\]%

\break

\noindent
then he proposed the following inequalities as a conjecture:%
\[
\left( \frac{\arcsin x}{x}\right)^{\!2}
\!+
\frac{\arctan x}{x}
<
2+\frac{\pi ^{2}+\pi -8}{\pi }x^{3}\arctan x,
\;\;\;
x \in \left(0,1\right)
\]%
and
\[
2\left(
\frac{\arcsin x}{x}
\right)
+
\frac{\arctan x}{x}
<
3 + \frac{5\,\pi-12}{\pi} \, x^3 \arctan x,
\;\;\;
x \in \left(0,1\right).
\]
Very recently, Male\v sevi\' c et al.
\cite{Malesevic_Banjac_Jovovic_2015} solved this open problem using
the same procedure - the natural approach method - associated to
upwards and downwards approximations of the inverse trigonometric
functions.

\medskip
Finally, we present other ways for approximating the function $%
\cos ^{2n}{\!x}$, $n\in \NN$. It is well known that the power series
of the function $\cos ^{2n}{\!x}$ converges to the function
everywhere on  $ \RR$. The power series of the function $\cos
^{2n}{\!x}$ is an alternating sign  series.  For example, for $n=1 $
and $x\in
\mathbb{R}
$, we have:
\[
\cos ^{2}{\!x}=1-x^{2}+\frac{1}{3}x^{4}-\frac{2}{45}x^{6}+\ldots =1+%
\displaystyle\sum_{k=0}^{\infty }{\displaystyle\frac{2^{2k-1}(-1)^{k}}{(2k)!}x^{2k}}.
\]%
Therefore, for the above power (Taylor) series it is not hard to
determine (depending on  $m$) which partial sums (i.e. Taylor
polynomials) $\displaystyle T_{m}^{\,\cos ^{2}\!x,\mbox{\scriptsize
$0$}}(x)$ become good downward or upward approximations of the function $%
\cos ^{2}{\!x}$ in a given interval \mbox{${\cal I}$}. Assuming the
following representation  of the function $\cos ^{2n}{\!x}$ in power
(Taylor) series
\[
\cos ^{2n}{\!x}%
=a_{0}^{(2n)}-a_{2}^{(2n)}x^{2}+a_{4}^{(2n)}x^{4}-a_{6}^{(2n)}x^{6}+\ldots ,
\]%
with  $a_{j}^{(2n)}>0$ $\left( j=0,2,4,6,\ldots \right) $, the power
(Taylor) series of function $\cos ^{2n+2}{\!x}$ will  be an
alternating sign series as follows:
\[
\begin{array}{rcl}
\cos ^{2n+2}{\!x}\! & \!=\! & \!\cos ^{2}{\!x}\cdot \cos ^{2n}{\!x} \\[1.5ex]
\! & \!=\! & \!\;\;\;\underbrace{\mathop{a_{0}^{(2n)}}}\limits%
_{a_{0}^{(2n+2)}} \\[3ex]
\! & \!\! & \!-\underbrace{\mathop{{\big (}a_{0}^{(2n)}+a_{2}^{(2n)}{\big )}}%
}\limits_{a_{2}^{(2n+2)}}x^{2} \\[3ex]
\! & \!\! & \!+\underbrace{%
\mathop{{\Big (}\mbox{\small
$\displaystyle\frac{1}{3}$}a_{0}^{(2n)}+a_{2}^{(2n)}+a_{4}^{(2n)}{\Big )}}}\limits%
_{a_{4}^{(2n+2)}}x^{4} \\[2ex]
\! & \!\! & \!-\underbrace{%
\mathop{{\Big (} \mbox{\small
$\displaystyle\frac{2}{45}$}a_{0}^{(2n)} + \mbox{\small
$\displaystyle\frac{1}{3}$}a_{2}^{(2n)} + a_{4}^{(2n)} + a_{6}^{(2n)}{\Big )}}}\limits%
_{a_{6}^{(2n+2)}}x^{6} \\[2ex]
\! & \!\! & \!+\,\ldots .
\end{array}%
\]%
\noindent with $a_{j}^{(2n+2)}>0$ ${\big (}j=0,2,4,6,\ldots {%
\big )}$.

\break

\medskip \noindent Therefore, in general, for the function $\cos^{2n}{\!x}$ it
is possible to determine, depending on the form of the real natural number $m
$, the upward (downward) Taylor approximations $\overline{T}^{\,\cos^{2n}{\!x%
},\mbox{\scriptsize $0$}}_{m}(x)$ {\Big (}$\underline{T}^{\,\cos^{2n}\!x,%
\mbox{\scriptsize $0$}}_{m}(x)${\Big )} that are all above (below) the
considered function in a given interval $\mathcal{I}$. Such estimation of
the function $\cos^{2n}{\!x}$ and the use of corresponding Taylor
approximations will be the object of future research.

\bigskip \noindent \textbf{Acknowledgements.} The first and the second
author was supported in part by the Serbian Ministry of Education,
Science and Technological Development, Projects ON 174032, III 44006
and TR 32023. The third author was  supported by a Grant of the
Romanian National Authority for Scientific Research, CNCS-UEFISCDI,
with the Project Number PN-II-ID-PCE-2011-3-0087.

\end{document}